\documentclass{elsart}
\input epsf
\usepackage{psfrag}

\newcommand{\bd}{\partial}

\newcommand{\QED}{\hfill $\sqcap \hspace{-3mm} \sqcup$}
\begin{document}
\begin{frontmatter}
\title{Visually Building Smale flows in $\bf S^3$}

\author{Michael C.~Sullivan}

\address{
Department of Mathematics (4408),
Southern Illinois University,
Carbondale, IL 62901, USA,
{\tt msulliva@math.siu.edu},
{\tt http://www.math.siu.edu/sullivan}}

\date{\today}

\begin{abstract}
A Smale flow is a structurally stable flow with one dimensional 
invariant sets. We use information from homology and {\em template\/} 
theory to construct, visualize and in some cases, classify, 
nonsingular Smale flows in the 3-sphere.
\end{abstract}


\end{frontmatter}

%%%%%%%%%%%%%%%%%%%%%%%%%%%%%%%%%%%
%                                 %
%           SECTION 0             %
%                                 %
%%%%%%%%%%%%%%%%%%%%%%%%%%%%%%%%%%%
\setcounter{section}{-1}
\section{Introduction} \label{sec0}

The periodic orbits of a flow in $S^3$ form knots. 
For {\em Morse-Smale flows\/} there are only finitely
many periodic orbits. Wada \cite{W} has classified 
all links that can be realized as a collection of
closed orbits of nonsingular Morse-Smale flows on $S^3$.
Further, Wada's scheme includes an indexing of the 
components of the link according to whether the 
orbit is an attractor, a repeller, or a saddle.

In a {\em Smale flow}, by contrast, the saddle sets
may contain infinitely many closed orbits, while the 
attractors and repellers must still be collections of
finitely many orbits. Franks \cite{Fra85} devised an
abstract classification scheme for nonsingular
Smale flows on $S^3$ using a device he called the
{\em Lyapunov graph\/} of a flow. Each vertex of a Lyapunov graph 
corresponds to an attractor, repeller or saddle set (the 
{\em basic sets\/} of the flow). The saddle vertices are labeled 
with an {\em incidence matrix\/} 
(determined non-uniquely by the first return map on a cross section).
A simple algorithm is used to decide if a given Lyapunov graph can 
be realized by a nonsingular Smale flow on $S^3$. 
However, the Lyapunov graph contains no explicit 
information about the embedding of the basic sets.
In contrast to Wada's study of Morse-Smale flows, Franks' work does
does not allow one to {\em see\/} Smale flows. 
It was our curiosity to visualize Smale flows that motivated
this paper. It is however worth noting that Wada's results
have provided tools for understanding bifurcations between
Morse-Smale flows \cite{CMAV} and it is likely that some of
our results may sheld light on bifurcations between Smale
flows and form Morse-Smale flows to Smale flows. Also see
\cite{GY} for an example.

The project of this paper is to visually construct 
examples of Smale flows and in some special contexts classify
all the possible embedding types. Our primary tools 
will be the theory of {\em templates}, branched 2-manifolds
which model the saddle sets \cite{BW2} and certain earlier 
results of Franks that do give some information about 
the embedding of closed orbits. Specifically,  computations
of linking numbers and Alexander polynomials are employed.

Sections \ref{sec1} and \ref{sec2} contain background information. 
Our main classification theorem (Theorem~\ref{thm_LS}) is in 
section \ref{sec3}. Various generalizations and applications follow
in sections \ref{sec4} and \ref{sec5}.
The author wishes to thank John Franks and Masahico Saito
for helpful conversations.

%%%%%%%%%%%%%%%%%%%%%%%%%%%%%%%%%%%
%                                 %
%           SECTION 1             %
%                                 %
%%%%%%%%%%%%%%%%%%%%%%%%%%%%%%%%%%%
\section{Knots and links} \label{sec1}

A knot $k$ is an embedding of $S^1$ into $S^3$. 
It is traditional to use $k$ to denote both the 
embedding function and the image in $S^3$.
A knot may be given an orientation. We will always
use a flow to induce an orientation on our knots.
The {\em knot group}
of $k$ is the fundamental group of $S^3 \backslash k$.
A link of $n$ components is an embedding of $n$
disjoint copies of $S^1$ into $S^3$.

Two knots $k_1$ and $k_2$ (or two links) are equivalent if there is an
isotopy of $S^3$  that takes $k_1$ to $k_2$. When we talk about a knot
we almost always mean its equivalence class, or {\em knot type}. 
A {\em knot diagram\/} is a projection of a knot or link into a plane 
such that any crossings are transverse. If the knot has been given an 
orientation the crossings are then labeled as positive or negative
according to whether they are left-handed or right-handed respectively.

The knot group can be calculated from a diagram and, unlike the diagram, 
is invariant. If a knot has a diagram with no crossings then it is called 
an {\em unknot\/} or a {\em trivial knot\/}. The following proposition 
will be of use to us.

\begin{prop}{\bf (The Unknotting Theorem, \cite[page 103]{R})} 
\label{prop_unknotting}
The knot group of $k$ is infinite cyclic if and only if $k$ is the 
unknot.
\end{prop}

\begin{figure}[htb]
	\begin{center}  \mbox{
		\epsfxsize=1in
		\epsfbox{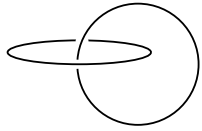}  
		}
	\end{center}
	\caption{The Hopf link}
	\label{fig_hopf}
\end{figure}

Given a diagram of a two component link $k_1 \cup k_2$ the 
{\em linking number\/} of
$k_1$ with $k_2$ is half the sum of the signs of each crossing of $k_1$
under $k_2$ and is denoted $lk(k_1, k_2)$.
The linking number is a link invariant, and thus is independent
of the choice of the diagram. The Hopf link shown in 
Figure~\ref{fig_hopf}. 
Its linking number is $\pm 1$, depending on the choice of orientations.

% @@@@@@@@@@@@@@@@@@@@@@@@@@@@ FIGURE @@@@@@@@@@@@@@@@@@@@@@@@@@@@@@
% @@@@@@@@@@@@@@@@@@@@@@@@@@@@@@@@@@@@@@@@@@@@@@@@@@@@@@@@@@@@@@@@@@@

\begin{figure}[htb]
	\begin{center}  \mbox{
		\psfrag{Prime}{Prime}
		\psfrag{Left-hand Trefoil}{}
		\psfrag{Right-hand Trefoil}{}
		\psfrag{2}{}
		\psfrag{S}{$S^2$}
		\psfrag{Square Knot}{}
		\psfrag{Composite}{Composite}
		\psfrag{p}{$p$}
		\psfrag{q}{$q$}
		\psfrag{g}{\hspace{-.04in}$\gamma$}
		\epsfxsize=4.5in
%		%\epsfysize=5in
		\epsfbox{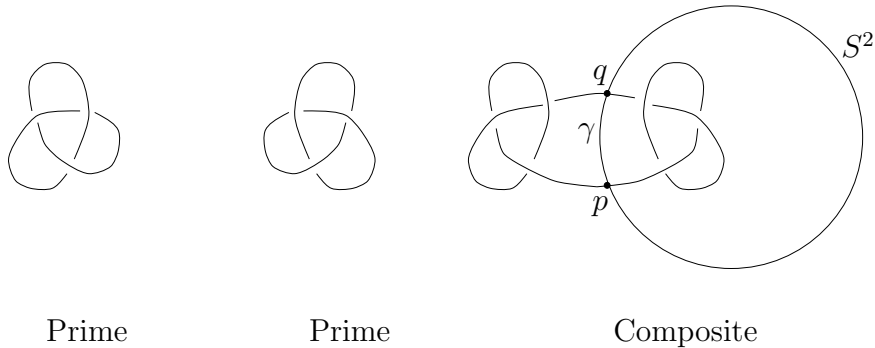}  
		}
	\end{center}
	\caption{The square knot is the sum of two trefoils.}
	\label{prime}
\end{figure}

A knot $k \subset S^3$ is {\em composite\/} if there exists a smooth
2-sphere $S^2$ such that $S^2 \cap k$ is just two points $p$ and $q$,
and if $\gamma$ is any arc on $S^2$ joining $p$ to $q$ then the knots
\[ k_1 = \gamma \cup ( k \cap \mbox{ outside of } S^2) \mbox{ and } \]
\[ k_2 = \gamma \cup ( k \cap \mbox{ inside of } S^2),   \] 
are each nontrivial, (i.e. not the unknot). We call $k_1$ and $k_2$
factors of $k$  and write $k = k_1 \# k_2.$
Of course the designation of the two components of $S^3/S^2$ as 
``inside'' and ``outside'' is arbitrary. This implies 
$k_1 \# k_2 = k_2 \# k_1$.
We call $k$ the connected sum of $k_1$ and $k_2$.
If a nontrivial knot is not composite then it is {\em prime}.
Figure~\ref{prime} gives an example. It shows how to factor the {\em
square knot\/} into two trefoils. Trefoils are prime. It was shown by
Schubert \cite[Chapter 5]{BZ} that any knot can be factored 
uniquely into primes, up to order. The unknot serves as a unit.

% @@@@@@@@@@@@@@@@@@@@@@@@@@@@ FIGURE @@@@@@@@@@@@@@@@@@@@@@@@@@@@@@
% @@@@@@@@@@@@@@@@@@@@@@@@@@@@@@@@@@@@@@@@@@@@@@@@@@@@@@@@@@@@@@@@@@@

A knot which fits on a torus and wraps about it $p$ times longitudinally
and $q$ times meridianally ($p$ and $q$ must be relatively prime), 
is called a $(p,q)$ {\em torus knot}.
If $k$ is any knot, then a $(p,q)$ {\em cable\/} about $k$, where $p$ 
and $q$
are relatively prime, is defined as follows.
Let $N(k)$ be a solid torus neighborhood of $k$, whose core is $k$.
Let $l$ be a {\em standard longitude\/} of $\partial N(k)$ for $k$,
i.e. $lk(l,k)=0$. 
Now consider a torus $T$ with a $(p,q)$ knot on it. 
Let $h:T \rightarrow \partial N$ be a homeomorphism that takes a 
standard longitude of $T$ to $l$. The image of $(p,q)$ under this 
map is said to be a $(p,q)$ cable of $k$ or, $(p,q)k$.

We shall extend the usual cabling notation a bit. Let
$(0,1)k$ be a meridian of $k$, $(1,0)k$ be a standard 
longitude of $k$, and $(0,0)k$ be a loop bounding a
disk on $T$. (A curve on a surface is {\em inessential\/} 
if it bounds a disk in the surface, and is {\em essential\/} 
if it does not.) 

Thus, we have tools with which to build up new knots from old ones.
In Wada's paper and in Sections 4 and 5 here, one builds new flows
from old ones using processes based in part on taking connected
sums and forming cables.

The cabling construction has been generalized in two ways. A knot $k'$ 
is a {\em satellite\/} of a given knot $k$ if $k'$ lives inside a 
tubular neighborhood $k$ and meets every meridianal disk. 
If the orientation of $k'$ is always roughly the same as that of $k$
(i.e. there is a fibration of $N(k)$ by meridianal disks which are
always transverse to $k'$), we say that $k'$ is a {\em generalized 
cable} of $k$.
It is known that the satellite of a nontrivial knot is nontrivial
\cite{R}, a fact which we shall make use of in Case 3 in the 
proof of Theorem~\ref{thm_LS}.

Knot polynomials form an important class of knot and link
invariants. We shall make use of the first known knot 
polynomial, the Alexander polynomial. It can be readily 
calculated from either a knot diagram or the fundamental
group. The latter approach will be of special importance
to us. The reader who wishes
to check our polynomial calculations should be able to find all 
he or she needs in \cite{FC}.

%%%%%%%%%%%%%%%%%%%%%%%%%%%%%%%%%%%
%                                 %
%           SECTION 2             %
%                                 %
%%%%%%%%%%%%%%%%%%%%%%%%%%%%%%%%%%%
\section{Dynamics of flows} \label{sec2}

A $C^1$ flow $\phi_{t}$ on a compact manifold $M$ is called {\em
structurally stable\/} if any sufficiently close approximation 
$\psi_{t}$ in the $C^1$ topology is {\em topologically equivalent}, 
that is if there exists a homeomorphism $h:M \rightarrow M$ taking 
orbits of $\phi_{t}$ to orbits of $\psi_{t}$, preserving the flow 
direction. Structurally stable $C^1$ flows have a hyperbolic
structure on their chain-recurrent sets \cite{Hay97}. 
We define these concepts next.

A point $x \in M$ is {\em  chain-recurrent\/} for $\phi_{t}$  if for
every $\epsilon > 0$ and $T>0$ there exists a chain of points $x =
x_{0}, \dots,  x_{n} = x$ in $M$, and real numbers $t_{0}, \dots,
t_{n-1}$ all bigger than $T$ such that $d(\phi_{t_{i}}(x_{i}),x_{i+1})
< \epsilon$ when ever $0 \leq i \leq n-1$.  The set of all such points
is called the chain-recurrent set $\mathcal R$.  It is a compact set
invariant under the flow.

A compact invariant set $K$ for a flow $\phi_{t}$  has a {\em
hyperbolic structure\/} if the tangent bundle of $K$ is the 
Whitney sum of three bundles $E^{s}$, $E^{u}$, and $E^{c}$ each of which
invariant under $D\phi_{t}$ for all $t$. Furthermore, the vector field
tangent to $\phi_{t}$ spans $E^{c}$ and there exist real numbers 
$C>0$ and $\alpha > 0$
such that
\[ \| D\phi_{t}(v) \| \leq C e^{-\alpha t} \| v \| 
\mbox{   for $t \geq 0$ and $v \in E^{s}$},	\]
\[ \| D\phi_{t}(v) \| \leq C e^{\alpha t} \| v \| 
\mbox{   for $t \leq 0$ and $v \in E^{u}$}.	\]

We also define the local stable and unstable manifolds
associated to an orbit $\bf O$. They are respectively,
\[ W^{s}_{\mbox{\tiny loc}}({\bf O}) = \bigcup_{x \in \bf O} \{ y \in
M | d(\phi_t (x),\phi_t (y)) \rightarrow 0 \mbox{ as } t \rightarrow 
\infty
\mbox{ and } d(x,y) \leq \epsilon \} \]
and
\[ W^{u}_{\mbox{\tiny loc}}({\bf O}) = \bigcup_{x \in \bf O} \{ y \in
M | d(\phi_t (x),\phi_t (y)) \rightarrow 0 \mbox{ as } t \rightarrow 
- \infty 
\mbox{ and } d(x,y) \leq \epsilon \}. \]
The global stable and unstable manifolds are defined similarly by
removing the condition that $d(x,y) \leq \epsilon$.

It was shown by Smale that if the chain-recurrent set $\mathcal R$ 
of flow has a hyperbolic structure then $\mathcal R$
is the union of a finite collection of disjoint invariant compact 
sets called the {\em basic sets}. Each basic set $\mathcal B$ 
contains an orbit whose closure contains $\mathcal B$. The periodic 
orbits of a basic set $\mathcal B$ are known to be dense in 
$\mathcal B$. A basic may either consist of a single closed orbit or it
may contain infinitely many closed orbits and infinitely 
other nonperiodic chain-recurrent orbits. In the later case any
cross section is a Cantor set and the first return map is a
{\em subshift of finite type}. In the former case any cross section
is a finite number of points but the first return map is still a
(trivial) subshift of finite type. Thus, each basic set is a 
suspension of a subshift of finite type. A nontrivial basic 
will be called {\em chaotic}.

\begin{defn}
A flow $\phi_t$ on a manifold $M$ is called a {\em Smale flow\/} 
provided
\begin{itemize}
\item [(a)] the chain-recurrent set $\mathcal R$ of $\phi_t$ has a 
hyperbolic structure,
\item [(b)] the basic sets of $\mathcal R$ are one-dimensional, and
\item [(c)] the stable manifold of any orbit in $\mathcal R$ has 
transversal intersection with the unstable manifold of any other 
orbit of $\mathcal R$.
\end{itemize}
\end{defn}

Most references allow for fixed points but we will be working primarily
with nonsingular flows.
Smale flows on compact manifolds are structurely stable under $C^1$
perturbations but are not dense in the space of $C^1$ flows. It is 
easy to see that for dim $M$ = 3, each attracting and repelling
basic set is a closed orbit. The admissible 
saddle sets, however, may be chaotic. A Smale flow with no 
chaotic saddle sets is called a {\em Morse-Smale\/} flow.

For a chaotic saddle set of a Smale flow in a 3-manifold one can
construct a neighborhood that is foliated by local stable manifolds of
orbits in the flow. Collapsing in the stable direction produces a
branched 2-manifold. With a semi-flow induced from the original flow,
this branched 2-manifold becomes what is known as a {\em template}. The
template models the basic saddle set in that the saddle set itself can
be recovered from the template via an inverse limit process and that
any knot or link of closed orbits in the flow is  smoothly isotopic 
to an equivalent knot or link of closed orbits in the template's
semi-flow. The proof of this is due to Birman and Williams \cite{BW2}
and can also be found in \cite[Theorem 2.2.4]{GHS}.

A key tool in the analysis of hyperbolic flows is the
concept of a Markov partition.
We refer the reader to 
\cite{Franksbook} for details. In our context a Markov partition
is a finite, disjoint collection of disks transverse to a
basic set of a flow. Each orbit of the basic set must
pass through some element of the Markov partition in forward time.

\begin{defn}
Given a Markov partition $\{m_1, \dots, m_n \}$ for
a suspended subshift of finite type 
with first return map $\rho$ 
we define the corresponding $n \times n$
{\em incidence matrix\/} $A$, by 
\[ A_{ij} = \left\{ \begin{array}{cl}
1	&	\mbox{if } \rho (m_i) \cap m_j \neq \phi   \\ 
0 	&	\mbox{otherwise.}
		\end{array} \right.	\]
\end{defn}

The incidence matrix, like a knot diagram, is not invariant but 
does contain invariant information. 
We can encode additional information about the embedding 
of a basic set by modifying the incidence matrix. 

\begin{defn} 
Given a Markov partition
for a basic set with first return map
$\rho$, assign an orientation to each partition element.
If the partition is fine enough the function 
\[ O (x) = \left\{ \begin{array}{cl}
	+1 & \mbox{if $\rho$ is orientation preserving at $x$}, \\
	-1 & \mbox{if $\rho$ is orientation reversing at $x$},
\end{array} \right. \]
is constant on each partition element. The {\em structure matrix\/} 
$S$ is then defined by $S_{ij} = O (x) A_{ij}$, where $x$ is any point
in the $i$-th partition element. 
\end{defn}

The next proposition was proved by Franks in \cite{Fra77}.

\begin{prop}[Franks, 1977] \label{prop_link}
Let $\phi_t$ be a Smale flow with a single attracting closed orbit $a$,
and a single repelling closed orbit $r$, with saddle sets 
$\Lambda_1, \dots, \Lambda_n$. 
Then the absolute value of the 
linking number of $a$ and $r$ is given by 
\[ 
 	|lk(a,r)| = \prod_{i=1}^n |\det (I-S_i)|, 
\]
where $S_1, \dots, S_n$ denote the respective structure matrices
of the saddle sets.
\end{prop}

The manner in which a saddle set ``links'' a closed orbit $k$, 
is described by modifying the structure matrix $S$ to 
form a {\em linking matrix\/} $L_k$. Consider a Markov partition, 
$\{m_1, \dots, m_n \}$, of the saddle with incidence matrix $A$.
Pick a base point $b$ in $S^3 - k$ and paths $p_i$ from
$b$ to $m_i$, also in $S^3 - k$. For each $a_{ij} \neq 0$ let 
$\gamma_{ij}$ be a segment of
the flow going from $m_i$ to $m_j$ without meeting any of the other
partition elements. Now form a loop consisting of $\gamma_{ij}$, 
$p_i$, $p_j$ and,
if needed, short segments in  $m_i$ and in $m_j$. If the disks
have been chosen small enough, then the linking number of any such 
loop with $k$ depends only on $i$ and $j$.
That this can always be done is shown in \cite{Fra81}.

\begin{defn}
The {\em linking matrix\/} $L$ associated with a suitable 
Markov partition for a given closed orbit $k$ is defined to by
\[ L_{ij} = S_{ij} t^q, \]
where $q$ is the linking number of the loops formed from segments 
connecting $m_i$ to $m_j$ and $k$ as described above.
\end{defn}

The following proposition is a special case of Theorem 4.1 in 
\cite{Fra81}.

\begin{prop} [Franks, 1981] \label{prop_alex}
Let $\phi_t$ be a Smale flow in $S^3$ with one attracting closed 
orbit $a$, one repelling closed orbit $r$, and a single saddle 
set $s$. Let $L_a$ and $L_r$ be the linking matrices for $s$ with 
respect to $a$ and $r$ respectively.
Then the Alexander polynomials of $a$ and $r$ are given by
$ \Delta_a(t) = \det(I-L_a)$ and $\Delta_r(t) = \det(I-L_r)$
respectively, up to multiples of $\pm t$.
\end{prop}

Finally, we record a proposition about Morse-Smale flows. 
It is an obvious corollary to Wada's theorem \cite{W}, though 
it is easy to prove directly.

\begin{prop} \label{prop_hopf}
In a (nonsingular) Morse-Smale flow on $S^3$ with exactly two closed 
orbits, the link of closed orbits forms a Hopf link with
one component a repeller and the other an attractor.
\end{prop}

%%%%%%%%%%%%%%%%%%%%%%%%%%%%%%%%%%%
%                                 %
%           SECTION 3             %
%                                 %
%%%%%%%%%%%%%%%%%%%%%%%%%%%%%%%%%%%
\section{Lorenz-Smale flows} \label{sec3}

\begin{figure}[ht]
        \begin{center}  \mbox{
                \epsfxsize=3in
                %\espfysize=3in
                \epsfbox{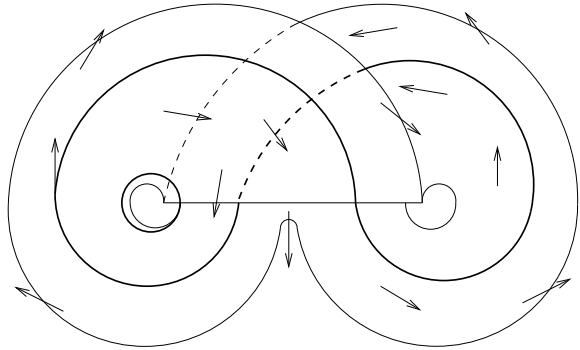}
                }
        \end{center}
        \caption{The Lorenz template}
        \label{Lorenz}
\end{figure}

By a {\em simple Smale flow} we shall mean a Smale flow with three
basic sets: a repelling orbit $r$, an attracting orbit $a$, 
and a nontrivial saddle set. In this section we will show how to 
construct simple Smale flows in which 
the saddle set $l$ can be modeled by an embedding of the
Lorenz template shown in Figure~\ref{Lorenz}. That is, there is an
isolating neighborhood of the saddle set $l$ foliated by local stable
manifolds of the flow, such that when we collapse out in the stable
direction we get an embedding of the Lorenz template. Call such flows 
{\em  Lorenz-Smale flows}. In this section we classify all
possible Lorenz-Smale flows.

% @@@@@@@@@@@@@@@@@@@@@@@@@@@@ FIGURE @@@@@@@@@@@@@@@@@@@@@@@@@@@@@@
% @@@@@@@@@@@@@@@@@@@@@@@@@@@@@@@@@@@@@@@@@@@@@@@@@@@@@@@@@@@@@@@@@@@

\begin{figure}[ht]
        \begin{center}  \mbox{
		\psfrag{Exit Set}{Exit Set}
		\psfrag{Entence Set}{Entrance Set}
		\psfrag{S}{\small $S$}
		\psfrag{S'}{\small $S'$}
		\psfrag{X}{\small $X$}
		\psfrag{X'}{\small $X'$}
		\psfrag{Y}{\small $Y$}
		\psfrag{Y'}{\small $Y'$}
                \epsfxsize=5in
                %\espfysize=3in
                \epsfbox{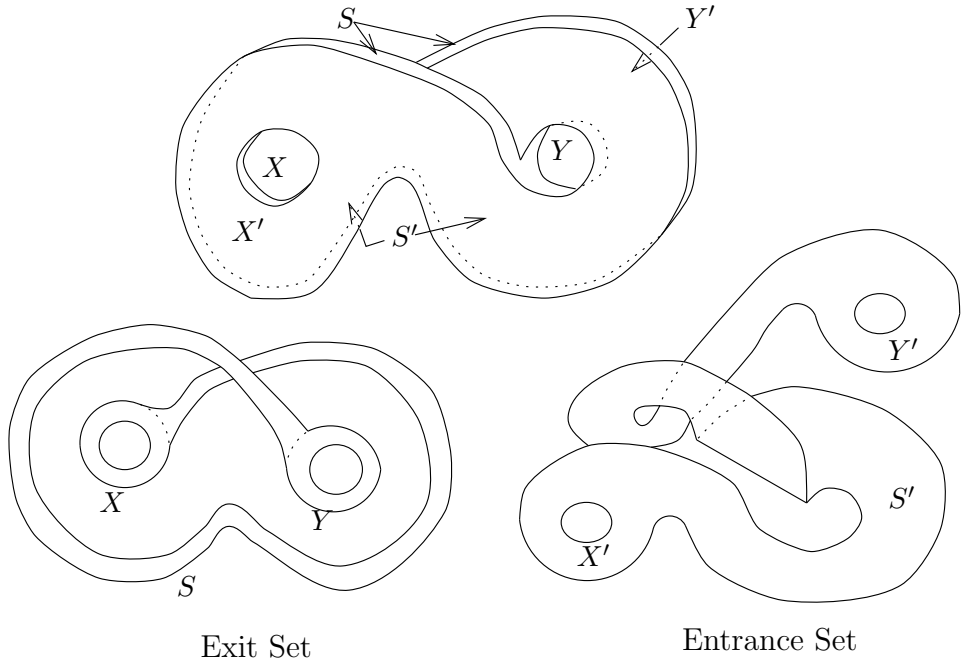}
                }
        \end{center}
        \caption{Neighborhood of a Lorenz saddle set (top) with the 
		exit set (bottom left) and entrance set (bottom right)}
        \label{Lnbhd}
\end{figure}

We shall call the isolating neighborhood of the saddle set $L$.
The set of points of $\bd L$ where the flow is transverse outward is 
called the {\em exit set\/}. The backward orbits of these points 
approach
orbits in the saddle set $l$. The exit set consists of two annuli, 
$X$ and $Y$, connected by a long strip $S$. The core of the exit 
set is homeomorphic to the boundary of the Lorenz template.
We shall refer to the cores of $X$ and $Y$, as $x$ and
$y$, respectively. See Figure~\ref{Lnbhd}. The set of points of $\bd L$ 
where the flow is transverse inward  is called the {\em entrance set\/}.
The entrance set also consists of two annuli connected by a long strip. 
Denote these two annuli by $X'$ and $Y'$ and their respective cores 
$x'$ and $y'$ and the connecting strip by $S'$. See Figure~\ref{Lnbhd}. 
The intersection of the closures of the entrance and exit sets consists 
of three closed curves where the flow is tangential to $\bd L$.
Although the entrance set is harder to visualize, its topological
type can be determined by an Euler characteristic argument. Notice 
that $x$ and $x'$ are isotopic to $\overline{X}\cap\overline{X'}$
and hence have the same knot type. Similarly, $y$ and $y'$ must have
the same knot type.

% @@@@@@@@@@@@@@@@@@@@@@@@@@@@ FIGURE @@@@@@@@@@@@@@@@@@@@@@@@@@@@@@
% @@@@@@@@@@@@@@@@@@@@@@@@@@@@@@@@@@@@@@@@@@@@@@@@@@@@@@@@@@@@@@@@@@@

Take tubular neighborhoods of $a$ and $r$ and denote them by $A$ and
$R$ respectively.  

To build a Smale flow from these building blocks, we first attach the 
closure of the exit set of $L$ to $\bd A$. This gives a vector 
field on a new 3-manifold pointing inward along its entire boundary. 
That this can be done smoothly was shown by Morgan in \cite{Mor78}. 
Next attach $\bd R$ to the boundary of $A \cup L$ so that the union 
$A \cup L \cup R$ is $S^3$.

We know from
Proposition~\ref{prop_link} that the linking number between $r$ and $a$
must be $\pm 1$, but what types of knots can $a$ and $r$ be? What are
all the different ways the saddle set can be embedded so as to still 
have the Lorenz template as a model? This later question is made 
precise by asking, what types of knots can $x$ and $y$ be? Can they 
be linked? Can the annuli $X$ and $Y$ have any number of twists or 
are there restrictions? (This last question is equivalent to finding 
the {\em framing\/} of a knot.)

To answer these questions we use the following framework. First, we
study what may happen when the attaching of the exit set of $L$ to
$\bd A$ is such that $x$ and $y$ are both inessential in $\bd A$.
Then we investigate the case where one of them, say $y$, is essential
but the other is still inessential. Finally we consider the case where
both $x$ and $y$ are essential in $\bd A$. The results of this
analysis are then consolidated into the statement of
Theorem~\ref{thm_LS}. Our classification scheme is only up to
isotopy of $S^3$, plus mirror images and flow reversal. Also, we
shall not be concerned with the orientation, i.e. the flow direction, 
of $a$ or $r$, since their orientations can be easily reversed by 
modifying the flow in a tubular neighborhoods of $a$ and $r$,
\cite[first paragraph]{W}.

\begin{thm} \label{thm_LS}
For a Lorenz-Smale flow in $S^3$  the following and only the
following
configurations are realizable. The link $a \cup r$ is either a Hopf
link or a trefoil and meridian. In the latter case the saddle set is
modeled by a standardly embedded Lorenz template, i.e. both bands are
unknotted, untwisted, and unlinked, with the core of each band a 
meridian
of the trefoil component of $a \cup r$. In the former case there are
three possibilities: (1) The saddle set is standardly embedded. (2)
One band is twisted with $n$ full-twists for any $n$, but remains
unknotted and unlinked to the other band, which must be unknotted and
untwisted. (3) One band is a $(p,q)$torus knot, for any pair of coprime
integers, with twist  $p+q-1$. The other band is unknotted, untwisted
and unlinked to the knotted one.
\end{thm}

\begin{pf}
The proof is divided into three cases.

{\bf CASE 1:} 
Suppose both $x$ and $y$ are inessential in $\bd A$.
It follows that $X$ and $Y$ are untwisted, that is the linking number
between each of the two components of $\bd X$ and of $\bd Y$ is zero.
It is also obvious that $x$ and $y$ are unknotted and unlinked.

There are two subcases to consider. It could be that $x$ and
$y$ are concentric in $\bd A$, or it could be that they are not. That
is $x$ and $y$ may or may not form the boundary of an annulus in 
$\bd A$.
In Figure~\ref{LonaBall} we show that both cases can be realized. 
The neighborhood $L$ is attached to a 3-ball $B$ along the closure of 
the exit set of $L$. Figure~\ref{LonaBall} also shows two ways one 
might 
attach handles to the 3-ball so as to turn it into a solid torus. 
Suppose we attach the handle to  the small disks marked $C$ and $C'$ 
in the manner shown. Call the resulting solid torus $A_1$. If we take 
$L \cup A_1$ the result is still a solid torus, and the complement in 
$S^3$ is just another solid torus, $R_1$. We can now build a Smale flow 
with an attractor in $A_1$, a repeller in $R_1$ and a Lorenz saddle set 
in $L$. Upon further inspection the reader should be able to see 
that $x$ and $y$ are concentric.

Now, instead on attaching a handle at $C$ and $C'$, we attach one
to $B$ and $B'$ as shown again in Figure~\ref{LonaBall}. This time
call the solid torus obtained $A_2$. As before $L \cup A_2$ is a solid
torus with solid torus complement in $S^3$. Thus, we have constructed 
another Lorenz-Smale flow with $x$ and $y$ inessential on a tubular 
neighborhood of the attractor.
Is it diffeomorphic to the Lorenz-Smale flow we constructed before?
To see that the answer is no, study the loops $x$ and $y$ again. They
are still both inessential, that is they both bound disks in $\partial
A_2$. But they are no longer concentric. This can be seen from careful
study of Figure~\ref{LonaBall}.

In both these examples the attractor and repeller form a Hopf link.
We claim that if both $x$ and $y$ are inessential in $\bd A$ then $a$
and $r$ must form a Hopf link. Since $R$ is a tubular neighborhood of 
a knot, the union of $A$ and $L$ is a knot complement. 
But we will show that $A \cup L$ is a solid
torus whose core has the same knot type as that of $a$. Thus, we could
remove $A \cup L$ from our flow and replace it with a solid torus
containing just an attractor and no saddle set. This gives us a 
nonsingular Morse-Smale flow on $S^3$ with just two closed orbits whose 
link type is the same as $a \cup r$. But by Proposition~\ref{prop_hopf} 
these must form a Hopf link.

We now show that $A \cup L$ is a solid torus of whose core is the same 
knot type as $a$. First, assume that $x$ and $y$ are not concentric. 
Then they bound disjoint disks $D_x$ and $D_y$ in $\bd A$. Expand $D_x$ 
and $D_y$, 
if needed, so that they contain all of $X$ and $Y$ respectively but remain
disjoint. Thicken $D_x$ and $D_y$ by pushing into $A$ a little, forming
two 3-balls $B_x$ and $B_y$, which are disjoint and do not meet the
orbit $a$. Let $L' = L \cup B_x \cup B_y$ and $A' = A - (B_x \cup
B_y)$. It is clear that $A'$ is a solid torus whose core is still $a$. 
The set $L'$ is a 3-ball and $L \cup A = L' \cup A'$. But the union
$L' \cup A'$ is taken along the disk  $S \cup \overline{(\bd B_x - D_x)}
\cup \overline{(\bd B_y - D_y)}$.  Thus, there is a deformation retract 
from $L'\cup A'$ to $A'$. This proves our claim for $x$ and $y$ not 
concentric.

Now suppose $x$ and $y$ are concentric and assume $x$ is inner
most. Then construct the 3-ball $B_x$ as before. Let $L' = L \cup B_x$
and $A' = A - B_x$. This time $L'$ is a solid torus. It is attached
to $A'$ along the annulus $Y \cup S \cup (B_x - D_x)$ whose core $y$
is inessential in $\bd A'$ and is a (1,0) longitude in $L'$. Thus we
can retract $L'$ to $Y$ and push it into $A'$ without changing the
knot type of $a$. This finishes the  proof of our claim.

% @@@@@@@@@@@@@@@@@@@@@@@@@@@@ FIGURE @@@@@@@@@@@@@@@@@@@@@@@@@@@@@@
\begin{figure}[ht]
        \begin{center}  \mbox{
                \epsfxsize=3in
                %\espfysize=3in
                \epsfbox{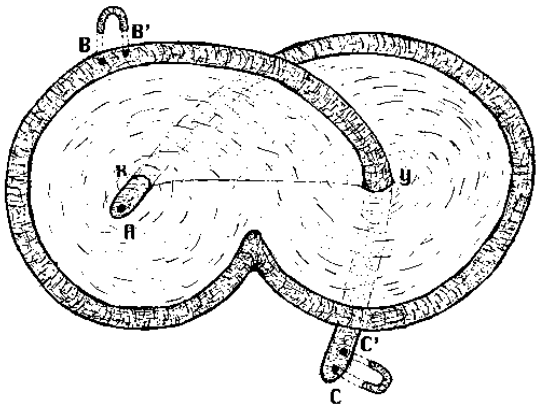}
                }
        \end{center}
        \caption{Lorenz saddle set neighborhood attached to a ball.}
        \label{LonaBall}
\end{figure}
% @@@@@@@@@@@@@@@@@@@@@@@@@@@@@@@@@@@@@@@@@@@@@@@@@@@@@@@@@@@@@@@@@@@

{\bf CASE 2:} 
Suppose that $x$ is inessential but that $y$ is essential.
The opposite case is similar. Again, it is clear that $X$ is untwisted
and that $x$ is unknotted and unlinked to $y$. We shall again show
that $a$ and $r$ must form a Hopf link. It then follows that since $y$ 
lives in a standardly embedded torus, $\bd A$, $y$ is a torus knot or 
unknot. If $y$ is a meridian (0,1), or a longitude (1,0) then
$Y$ is untwisted. If $y$ is an unknot $(1,q)$ or $(q,1)$ then $Y$ has
$q$ full twists. For nontrivial torus knots the twisting in $Y$ is
uniquely determined by the knot type of $y$. If $y$ is a $(p,q)$ torus
knot then the twist in $Y$ is $p+q-1$. 

That any $(p,q)$ torus knot 
can be realized by $y$ is shown by construction in Figure~\ref{Yess}.
One places a $(p,q)$ curve on a torus. Attach an annulus to this
curve along one boundary component. Add a ``Lorenz ear'' to form a
Lorenz template. Next a ``finger'' pushes out of the torus and snakes
along the boundary of the template and finally pokes through the $x$
loop. Thicken this complex up to get $A\cup L$. The repeller is then
placed as a meridian in the complement. An example with $y$ a $(2,1)$
curve is shown in Figure~\ref{YessEx}.

The argument that $a$ and $r$ must form a Hopf link is the same as in
the concentric subcase of Case 1 above. The core $y$ of the annulus
$Y$ is a $(1,q)$ cabling of the core of the solid torus $L'$. We can
foliate $L'$ with meridianal disks each of which meets $Y$ in an arc. 
Thus, $Y$ is a deformation retraction of $L'$. We then push 
$Y$ into $A'$.

% @@@@@@@@@@@@@@@@@@@@@@@@@@@@ FIGURE @@@@@@@@@@@@@@@@@@@@@@@@@@@@@@
\begin{figure}[htb] 
	\begin{center}  \mbox{
		\epsfxsize=3in
		\epsfbox{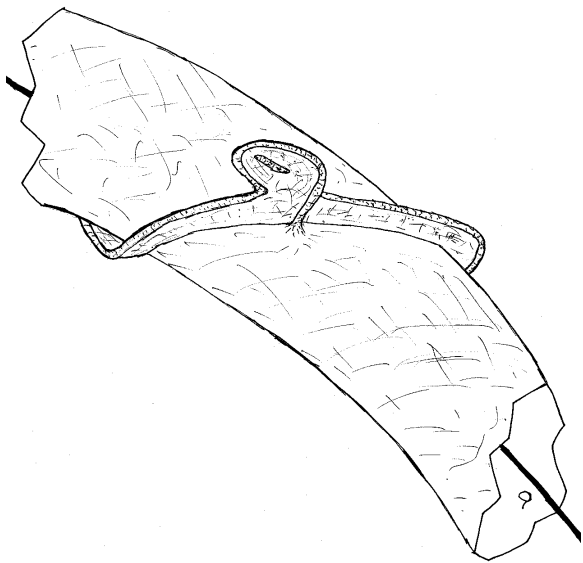}  
		}
	\end{center}
	\caption{$y$ is a $(p,q)$ cable of $a$}
	\label{Yess}
\end{figure}
% @@@@@@@@@@@@@@@@@@@@@@@@@@@@@@@@@@@@@@@@@@@@@@@@@@@@@@@@@@@@@@@@@@@

% @@@@@@@@@@@@@@@@@@@@@@@@@@@@ FIGURE @@@@@@@@@@@@@@@@@@@@@@@@@@@@@@
\begin{figure}[hbt]
        \begin{center}  \mbox{
                \epsfxsize=3in
                %\espfysize=3in
                \epsfbox{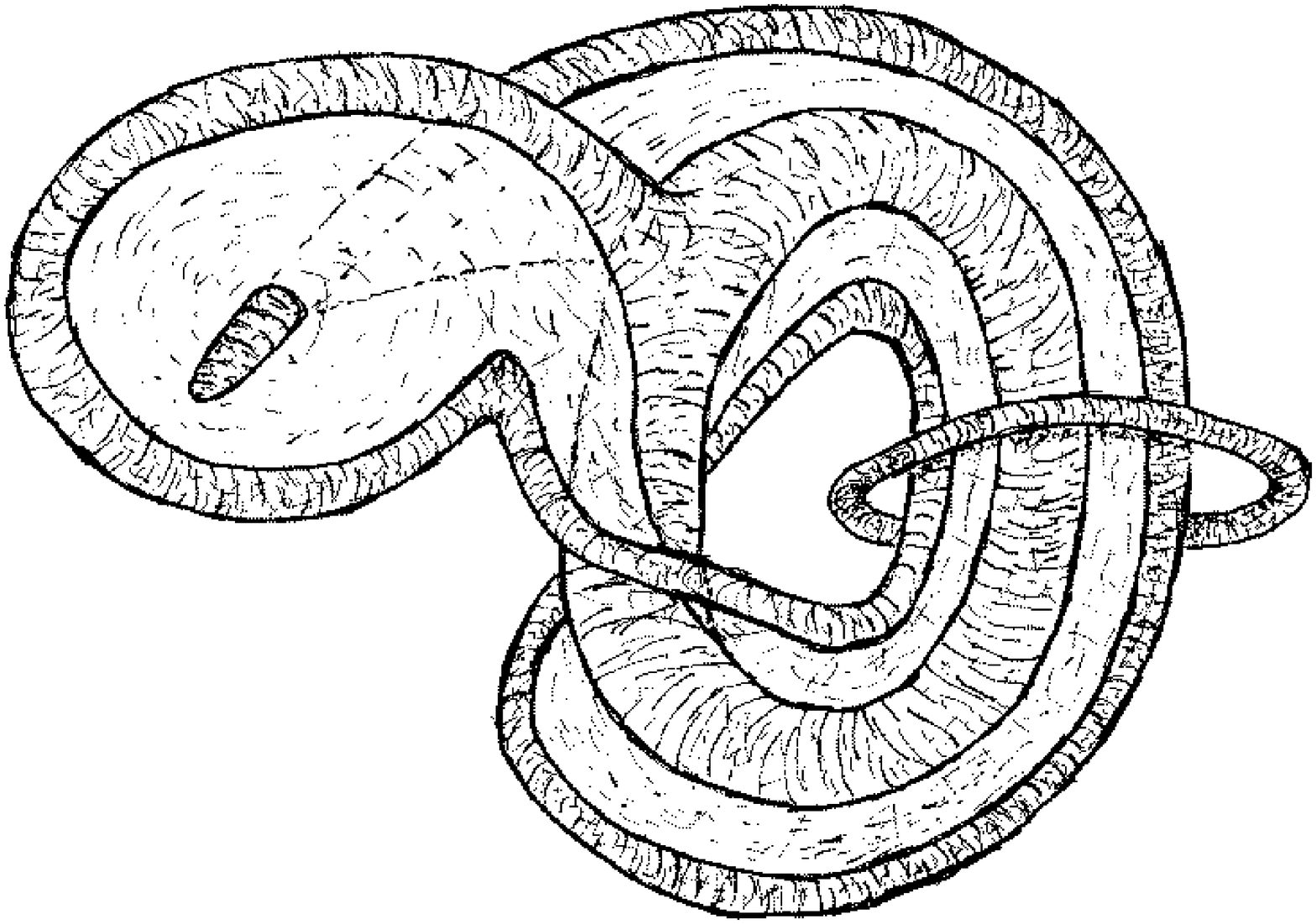}
                }
        \end{center}
	\vspace{.14in}
        \begin{center}  \mbox{
                \epsfxsize=3in
                %\espfysize=3in
                \epsfbox{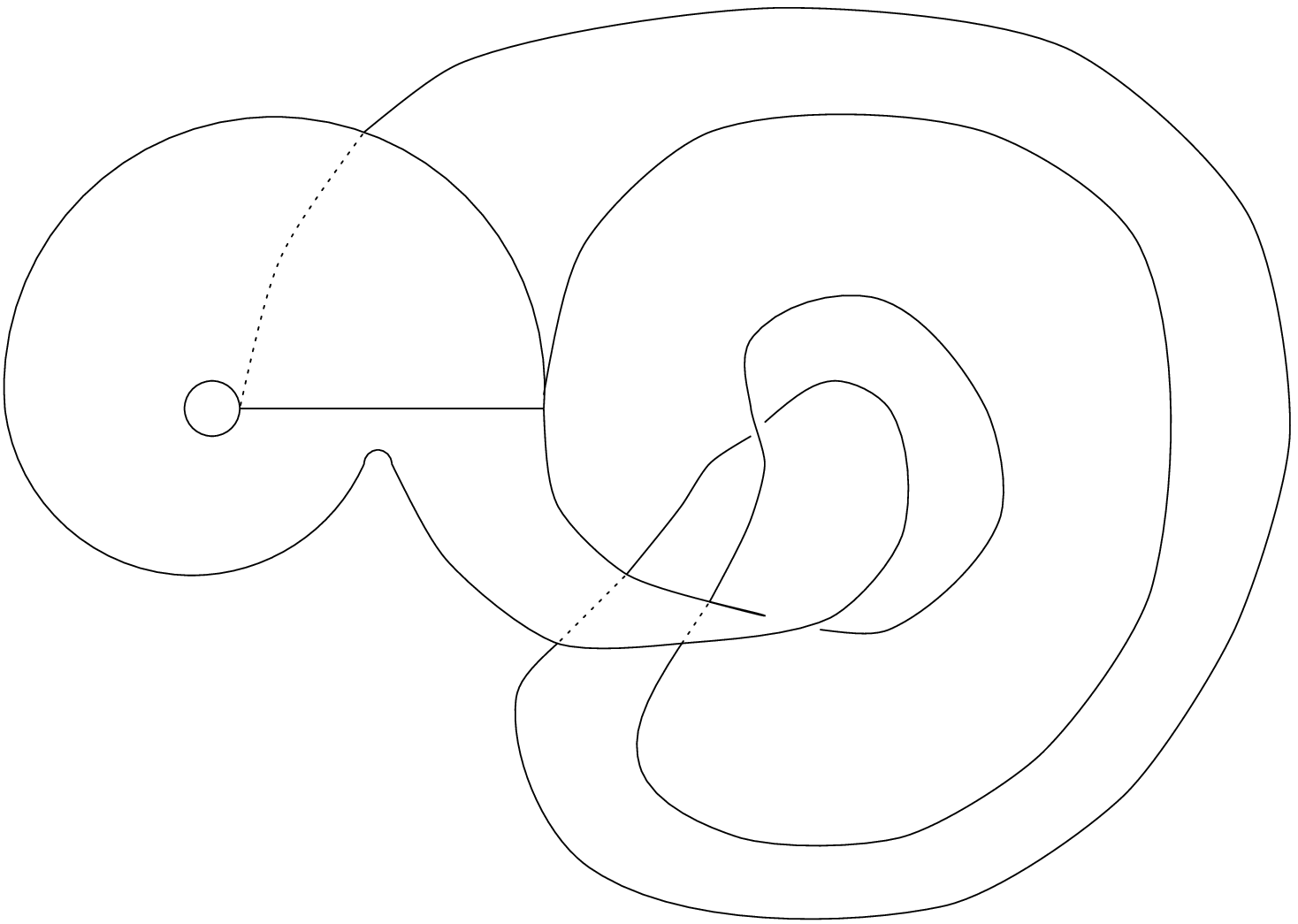}
                }
        \end{center}
        \caption{The top figure has an attractor in the fat tube
	and a repeller in the thin tube.  A template for the saddle
	set is shown below. The loop $y$ is a (2,1) torus curve. 
	Any $(p,q)$ torus curve, knotted or unknotted, can be realized.}
        \label{YessEx}
\end{figure}
% @@@@@@@@@@@@@@@@@@@@@@@@@@@@@@@@@@@@@@@@@@@@@@@@@@@@@@@@@@@@@@@@@@@

% @@@@@@@@@@@@@@@@@@@@@@@@@@@@ FIGURE @@@@@@@@@@@@@@@@@@@@@@@@@@@@@@
\begin{figure}[htb]
        \begin{center}  \mbox{
                \epsfxsize=3in
                %\espfysize=3in
                \epsfbox{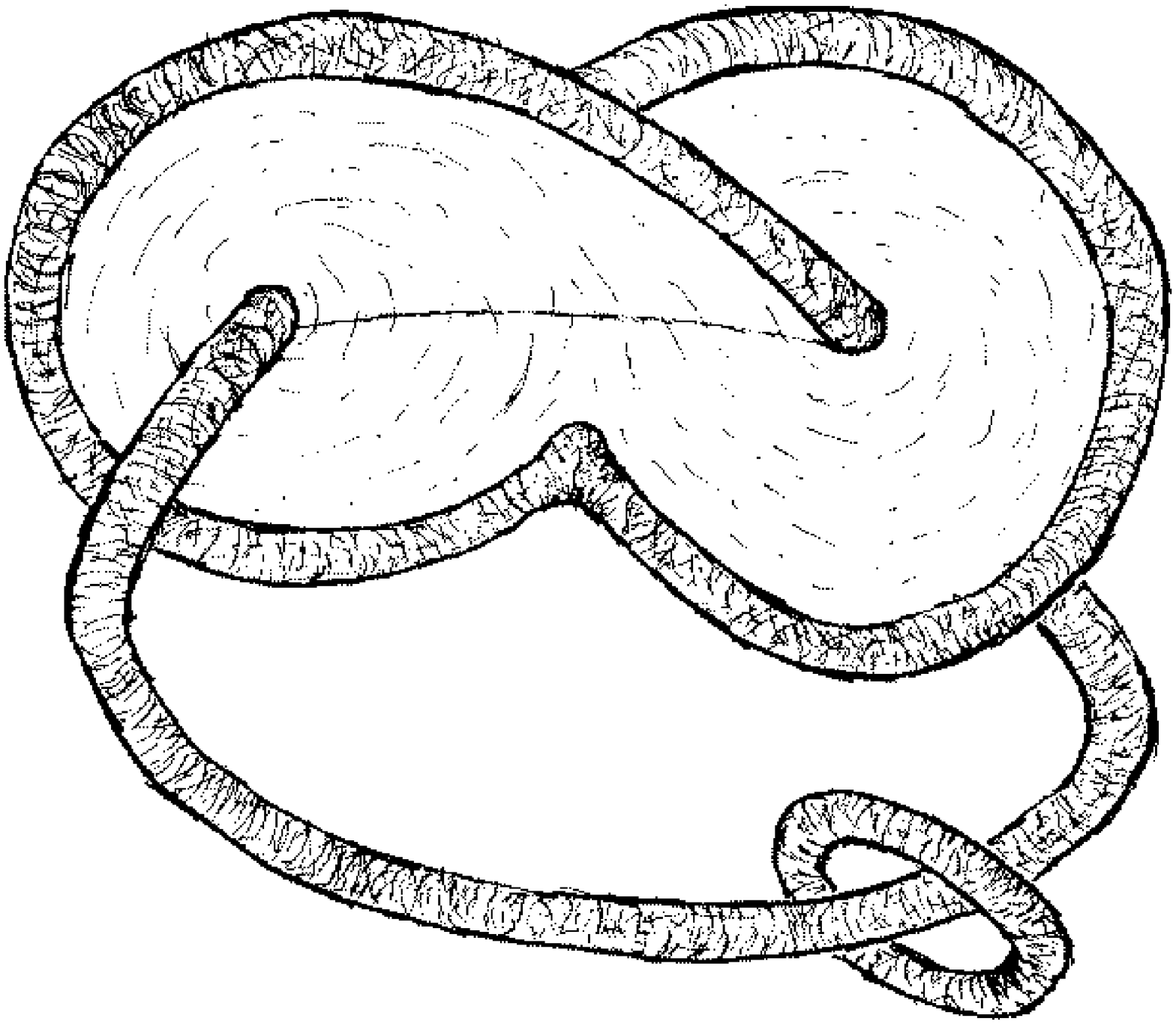}
                }
        \end{center}
        \caption{Lorenz-Smale flow with a trefoil attractor}
        \label{LandT}
\end{figure}
% @@@@@@@@@@@@@@@@@@@@@@@@@@@@@@@@@@@@@@@@@@@@@@@@@@@@@@@@@@@@@@@@@@@

{\bf CASE 3:}
 Suppose that both $x$ and $y$ are essential in $\bd A$.
In Figure~\ref{LandT} we give an example. The loops $x$ and $y$ are
meridians in $\bd A$. They are also
standard longitudes in $\bd R$. The attractor $a$ is a trefoil knot
while $r$ is unknotted and is a meridian of $a$. We claim that up to
mirror images and flow reversal, this is the only possible 
configuration.

Now consider the general setting. The attaching map from $\bd L$ to
$\bd A$ takes $x$ and $y$ to two copies of some $(p,q)$ cable
knot of the attractor, $a$. Here we allow $p$ or $q$ to be zero, but
not both.  Likewise, the attaching map from $\bd L$ to $\bd R$ takes
$x'$ and $y'$ to some pair of $(p',q')r$ knots. Of course $x$ and $y$
are respectively ambient isotopic to $x'$ and $y'$ within $\bd L$, so 
all four have the same knot type. 

For future reference, let $T = \left< x, y | xyx=yxy \right>$. A knot
$k$ with $\pi_1 (S^3 - k) = T$ is a left or right handed trefoil. See
\cite[Lemma 15.37, Corollary 15.23]{BZ}.

\begin{lem}
The Alexander polynomials of the attractor and repeller are
$\Delta_a=t^q-1+t^{-q}$ and $\Delta_r = t^{q'}-1+t^{-q'}$,
respectively
\end{lem}

\begin{pf}
If the linking matrix for $a$ is
\[ \left[ \begin{array}{cc}
	t^q	&	t^q	\\
	t^{-q}	&	t^{-q}
	  \end{array}	\right],	\]
then the result follows by Proposition~\ref{prop_alex}.
The only difficulty in determining the linking matrix is
the assignment of the signs to the powers of the $t$'s.
One can check our assignment explicitly for the $q=1$ case
by studying Figure~\ref{LandT}. In general, if the powers
are all of the same sign, the polynomial that results is
not symmetric in $t$, nor is any $\pm t$ multiple. But it is 
well known that the Alexander polynomial of a knot is
symmetric in $t$, up to multiples of $\pm t$. That is
$ \Delta(t) = \pm t^n \Delta(1/t), \mbox{ for some } n.$
See \cite{FC}. 
\QED \end{pf}

Not ready:

\begin{lem} \label{lem_11}
The fundamental groups of $L \cup A$ and $L \cup R$ are
$\left< x y | x^pyx^p=yx^py \right>$ and $\left< x y | 
x^{p'}yx^{p'}=yx^{p'}y \right>$, respectively. The Alexander
polynomials of $a$ and $r$ are $\Delta_a = t^{p'}-1+t^{-p'}$ and
$\Delta_r = t^p-1+t^{-p}$, respectively.
\end{lem}

\begin{pf}
We shall find $\pi_1(L\cup R)$ using the Seifert/Van Kampen Theorem.
The calculations of $\pi_1(L\cup A)$ are similar. The Alexander
polynomials can then be determined from Fox's Free Differential
Calculus \cite{FC}.

We must choose generators for $L$, $R$ and $L \cap R$. 
The generators for $L \cap R$ and $L$ are shown in
Figure~\ref{generators}. The base point $b,$ is in the ``middle'' 
of the strip $S'$. For $L$ we shall abuse notation slightly and call
the generators $x$ and $y$, as they are isotopic to the $x$ and $y$
loops, however, we do not use the orientation of the flow.
(By the proof of the previous lemma the images of $x$ and $y$ 
must wrap around $\bd R$ in the same direction.)
Denote the generators of $L \cap R$ by $w$ and $z$.
For $R$ we shall use a loop isotopic to $r$ but with base point 
$b \in \bd R \cap S'$.
Again we abuse notation and call this new loop $r$.

The fundamental groups of interest are then 
$  \pi_1 (R) = \left< r \right>,$ 
$  \pi_1 (L) = \left< x, y \right>,$ and 
$  \pi_1 (L \cap R) = \left< w, z \right>.$
The homomorphisms induced by inclusion maps are 
$\alpha: \pi_1 (L \cap R) \rightarrow \pi_1 (R)$  
and 
$\beta: \pi_1 (L \cap R) \rightarrow \pi_1 (L)$.
These give
$ \alpha (w) = r^p,$
$ \alpha (z) = \overline{r}^p,$
$ \beta (w) = yx\overline{y},$
and
$\beta (z) = xy\overline{x},$.
By Van Kampen's theorem 
$ \pi_1 (L \cup R) = \left< r, x, y | r^p = yx\overline{y}, \,
\overline{r}^p=x\overline{yx} \right>  \cong
\left<r, y | \overline{y}r^p\overline{y} 
= r^p\overline{y}r^p \right>$. 
\QED \end{pf}

% @@@@@@@@@@@@@@@@@@@@@@@@@@@@ FIGURE @@@@@@@@@@@@@@@@@@@@@@@@@@@@@@
\begin{figure}[hbt]
        \begin{center}  \mbox{
		\psfrag{x}{x}
		\psfrag{y}{y}
		\psfrag{w}{w}
		\psfrag{z}{z}
                \epsfxsize=3in
                %\espfysize=3in
                \epsfbox{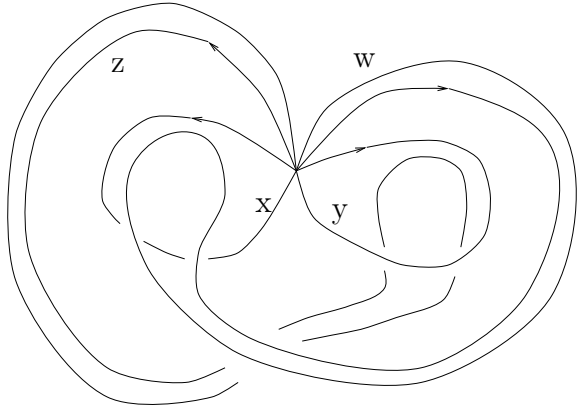}
                }
        \end{center}
        \caption{The generators of $\pi_1(R \cap L)$ and $\pi_1(L)$.}
        \label{generators}
\end{figure}
% @@@@@@@@@@@@@@@@@@@@@@@@@@@@@@@@@@@@@@@@@@@@@@@@@@@@@@@@@@@@@@@@@@@

It follows from Lemmas 10 and 11 that if $x$ is a $(p,q)$ curve on 
$\bd A$ then $x'$ is a $(\pm q,p)$ curve on $\bd R$.

If $p$ or $q$ is zero then the other is $\pm 1$ since the curve
is in a torus.
Now suppose $x$ is a $(0,\pm 1)$ curve on $\bd A$. Then $\pi_1(L \cup R)
\cong T$, and so $a$ is a trefoil knot. 
Since $x$ is a meridian of $a$ and $x$ is isotopic to $x'$ which
in turn is isotopic to $r$, we see that $r$ must be a meridian of $a$.

If $x$ is a $(\pm 1,0)$ curve on $\bd A$ then the rolls of $a$ and $r$
are switched.

\begin{figure}[htb]
	\begin{center}  \mbox{
		\epsfxsize=3in
		\epsfbox{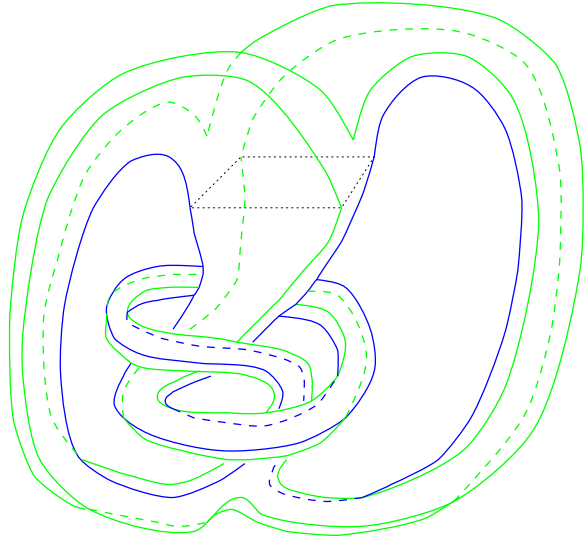}  
		}
	\end{center}
	\caption{The gray curve is $\gamma$.}
	\label{ktemp.eps}
\end{figure}

It is left only to show that $p$ and $q$ cannot both be nonzero.
It shall be useful to study the attaching of the exit set of $L$
to $\bd A$ in terms of the boundary curves of the exit set. 
They consist of three loops denoted as $\alpha$, $\beta$ and 
$\gamma$. We take $\alpha$ to be isotopic to $x$ and $\beta$ to
be isotopic to $y$. Then $\gamma$ is the remaining curve.
See Figure~\ref{ktemp.eps}. 
Our strategy is to show that if $p$ and $q$ are both nonzero then
$\gamma$ bounds a disk in $\bd A$ and that $\gamma$ is a nontrivial 
knot. This contradiction will then prove our claim. 

Let $\bd_+ L$ be the closure of the exit set of $L$. Clearly
$(\bd A \backslash \bd_+L) \cup \bd_+L$ is torus. Now $\alpha$
and $\beta$ bound an annulus $\alpha\beta$ in  
$(\bd A \backslash \bd_+L)$. Thus, $(\bd A \backslash \bd_+L)$ has
two components, the annulus $\alpha\beta$ and another component
we shall call $D$ which has a single boundary component $\gamma$.
Now, $\bd_+L \cup \alpha\beta$ is a torus with a disk removed.
Since $\bd_+L \cup \alpha\beta \cup D$ must be a torus, $D$ is a disk. 
Since this torus is embedded in $S^3$ it follows that $\gamma$, 
the boundary of $D$, is unknotted.

\begin{figure}[htb]
	\begin{center}  \mbox{
		\epsfxsize=4in
		\epsfbox{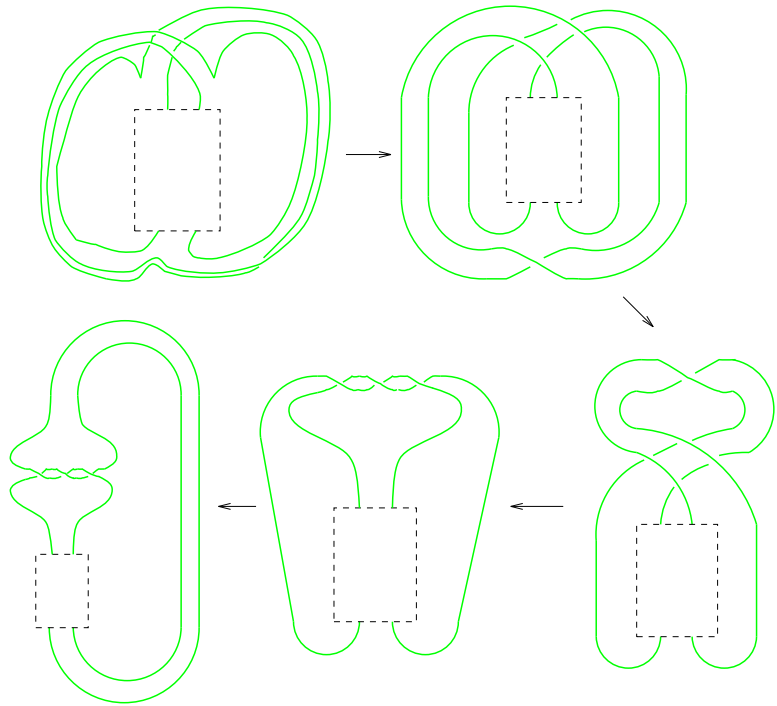}  
		}
	\end{center}
	\caption{$\gamma$ is a satellite of $\alpha$}
	\label{gamma}
\end{figure}

The Alexander polynomial calculations in Lemma 11 show that $a$ is knotted,
and thus even for the $(1,1)$ case $\alpha$ and $\beta$ are nontrivial knots.
Now since $\alpha$ and $\beta$ are parallel knots in $\bd A$ 
we can deform $L$ to appear as in Figure~\ref{ktemp.eps}.
By studying Figure~\ref{gamma} we see that $\gamma$ is a satellite
of $\alpha$. This implies $\gamma$ is nontrivial and completes the
proof of Theorem 9.
\QED \end{pf}

% @@@@@@@@@@@@@@@@@@@@@@@@@@@@ FIGURE @@@@@@@@@@@@@@@@@@@@@@@@@@@@@@
% @@@@@@@@@@@@@@@@@@@@@@@@@@@@@@@@@@@@@@@@@@@@@@@@@@@@@@@@@@@@@@@@@@@

% @@@@@@@@@@@@@@@@@@@@@@@@@@@@ FIGURE @@@@@@@@@@@@@@@@@@@@@@@@@@@@@@
% @@@@@@@@@@@@@@@@@@@@@@@@@@@@@@@@@@@@@@@@@@@@@@@@@@@@@@@@@@@@@@@@@@@

\begin{cor}
In any Lorenz-Smale flow there is a pair of unlinked saddle orbits.
\end{cor}

\begin{pf}
In the Lorenz template there are two fixed points in the first 
return map of the branch line. These correspond to a pair of 
unlinked unknots in the semi-flow if the Lorenz template is 
standardly embedded. In all of the embeddings
allowed for in Theorem~\ref{thm_LS} these two orbits remain unlinked.
Thus, by the Birman-Williams template theorem \cite{BW2}, the 
saddle set also contains a pair of unlinked orbits.
\QED \end{pf}

\begin{cor}
Consider a Smale flow on $S^3$ with an attracting fixed point, a
repelling fixed point, a Lorenz saddle set and no other basic sets.
There is only one possible configuration, and in it the template of
the saddle set is a standardly embedded Lorenz template. 
\end{cor}

\begin{pf}
 The proof is similar to Case 1 above.
\QED \end{pf}

%%%%%%%%%%%%%%%%%%%%%%%%%%%%%%%%%%%
%                                 %
%           SECTION 4             %
%                                 %
%%%%%%%%%%%%%%%%%%%%%%%%%%%%%%%%%%%
\section{Connected sums} \label{sec4}

Wada's classification theorem for Morse-Smale flows is based
on applying a series of {\em moves} to one or two existing 
Morse-Smale flows and building up new ones. Conclusion (a) of the
next theorem  establishes 
an operation that produces a new Smale flow from two existing
ones that is similar to Wada's move IV \cite{W}. 

\begin{thm} \label{thm_sums}
Let $\phi_1$ and $\phi_2$ be nonsingular Smale flows on $S^3$
such that (1) they each have only one attracting 
closed orbit with knot types $k_1$ and $k_2$ respectively, (2)
there is only one repelling closed orbit which is unknotted and 
is a meridian of the attractor and, (3) the repellers bound disks whose
interiors meet the chain-recurrent sets in a single point. 
Then (a) and (b) below hold true.

(a) There exists a nonsingular Smale flow on $S^3$ 
such that there is only one attracting closed orbit
which has knot type $k_1 \# k_2$, and there is only one
repelling closed orbit which is unknotted and is a meridian
of the attractor.

(b) There exists a nonsingular Smale flow on $S^3$
such that there is only one attracting closed orbit
which has knot type $k_1$, and there is only one repelling
closed orbit which has knot type $k_2$. Furthermore, the attractor
and the repeller have linking number one, and can be placed 
into solid tori whose cores are meridians of the each other.
\end{thm}

\begin{figure}[htb]
	\begin{center}  \mbox{
	\psfrag{D+1}{$D^+_1$}
	\psfrag{D-1}{$D^-_1$}
	\psfrag{D+2}{$D^+_2$}
	\psfrag{D-2}{$D^-_2$}
		\epsfxsize=4in
		\epsfbox{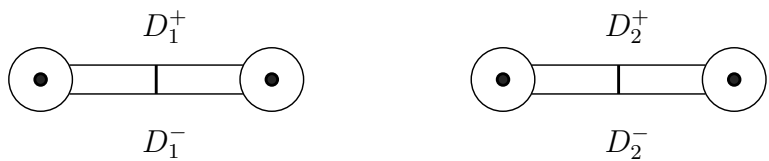}  
		}
	\end{center}
	\caption{Cut out these balls.}
	\label{fig_cut}
\end{figure}

\begin{figure}[htb]
	\begin{center}  \mbox{
		\epsfxsize=3in
		\epsfbox{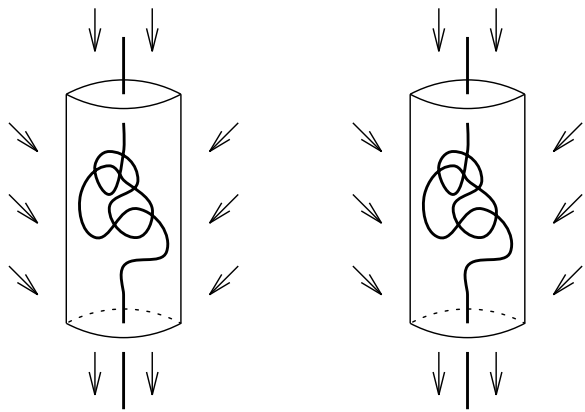}  
		}
	\end{center}
	\caption{Paste cylinders together}
	\label{fig_paste}
\end{figure}

\begin{pf}
The proofs are simple cut and paste arguments. Some details
are left to the reader.
For (a) let $V_i$, $i=1, 2$, be tubular neighborhoods of the 
repellers in $\phi_i$,  $i=1, 2$, respectively. Let $D_i$,
 $i=1, 2$, be the disks described in hypotheses (3). Thicken up
these disks a bit by taking cross product with a small interval
$I=[-1,1]$. We require that each $D_i \times I$ meet the 
corresponding $k_i$ in an unknotted arc.
Let $D_i^{\pm} = (D_i \times \{ \pm1 \}) \backslash V_i$. 
Choose the signs so that flows enter the thickened disks on 
the positive sides. 
See Figure~\ref{fig_cut}.
Delete from the 3-sphere of each flow the interior of the union 
of $V_i$ and $D_i \times I$, for the corresponding $i=1, 2$. 
We now have flows on two cylinders $C_i$,  $i=1, 2$. 
See Figure~\ref{fig_paste}. The boundary of $C_i$ is the union of 
$D_i^+$, $D_i^-$, and the annulus 
$A_i = \bd V_i \backslash D_i \times I$.
The flow exits $C_i$ on the interior of $D_i^+$, for $i=1, 2$.  
Glue $C_1$ to $C_2$ by identifying $D_1^+$ with $D_2^-$ and 
$D_2^+$ with $D_1^-$. This creates a solid torus, $V$. The flow 
induced on $V$ is inward on $\bd V = A_1 \cup A_2$. Further, 
the identifications can be chosen so that the flow on $V$ has an
attracting orbit with knot type $k_1 \# k_2$, assuming $V$ is 
standardly embedded. It is now clear how to construct the desired 
flow on $S^3$.

The proof of conclusion (b) is similar and in fact simpler and so
is left as an exercise. 
Note that hypotheses (3) is not required.
\QED \end{pf}

% @@@@@@@@@@@@@@@@@@@@@@@@@@@@ FIGURE @@@@@@@@@@@@@@@@@@@@@@@@@@@@@@
% @@@@@@@@@@@@@@@@@@@@@@@@@@@@@@@@@@@@@@@@@@@@@@@@@@@@@@@@@@@@@@@@@@@

% @@@@@@@@@@@@@@@@@@@@@@@@@@@@ FIGURE @@@@@@@@@@@@@@@@@@@@@@@@@@@@@@
% @@@@@@@@@@@@@@@@@@@@@@@@@@@@@@@@@@@@@@@@@@@@@@@@@@@@@@@@@@@@@@@@@@@

%%%%%%%%%%%%%%%%%%%%%%%%%%%%%%%%%%%
%                                 %
%           SECTION 5             %
%                                 %
%%%%%%%%%%%%%%%%%%%%%%%%%%%%%%%%%%%
\section{Attracting Links} \label{sec5}

The class of links which can arise in Morse-Smale flows
(nonsingular on $S^3$) is, according to Wada's Theorem,
quite limited. For Smale flows it is easy to construct
examples in which every knot and link can be  realized
simultaneously as a saddle orbits. This is a consequence of
the existence of {\em universal templates}, templates in which
all links all realized as closed orbits, \cite{Ghr97}; also see
\cite{GHS} and \cite{W98}.

Franks has shown that any link can be realized as an attractor
of a Smale flow \cite[Propsition 6.1]{Fra81}. 
In the proof 
the link is realized as a braid in an unknotted solid
torus whose entrance set is its entire boundary. 
Thus given a Smale flow with attractor $k$ we can
replace $k$ with any generalized cable of $k$, though a new saddle
set will typically be introduced. 

In \cite{GHS} it is shown that given a Smale flow $\phi$ with a saddle
set modeled by a template $T$ containing the closed orbit $k$,
there exists another Smale flow $\phi'$ with the same basic sets as
$\phi$ except that  $k$ is an attractor
and and the template $T$ has been replaced (as a model) with 
$T'$, a template formed by ``surgering'' $T$ along $k$ (a 
standard template operation).

Turning our attention to simple Smale flows we shall use a similar
idea to show that 
given any knot $k$ there exists a simple Smale flow with attractor 
$k$ and repeller a meridian of $k$. As a corollary of the construction 
we  can give a ``dynamics'' proof that Alexander polynomials 
multiply under connected sums.

\begin{thm} \label{thm_anyk}
For any knot $k$ there exists a simple Smale flow such that with 
attractor $k$ and repeller a meridian of $k$ that does not link any 
closed orbits in the saddle set.
\end{thm}

\begin{figure}[htb]
\psfrag{a}{a.}
\psfrag{b}{b.}
	\begin{center}  \mbox{
		\epsfxsize=3in
		\epsfbox{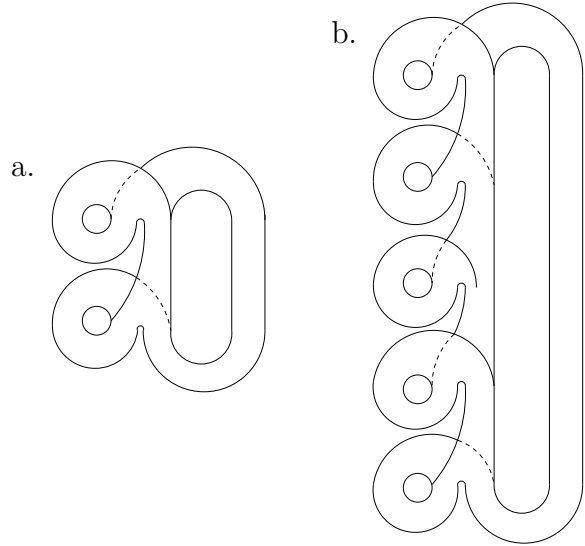}  
		}
	\end{center}
	\caption{Two templates containing all knots}
	\label{univ.eps}
\end{figure}

\begin{figure}[htb]
\psfrag{k}{$k$}
	\begin{center}  \mbox{
		\epsfysize=4in
		\epsfbox{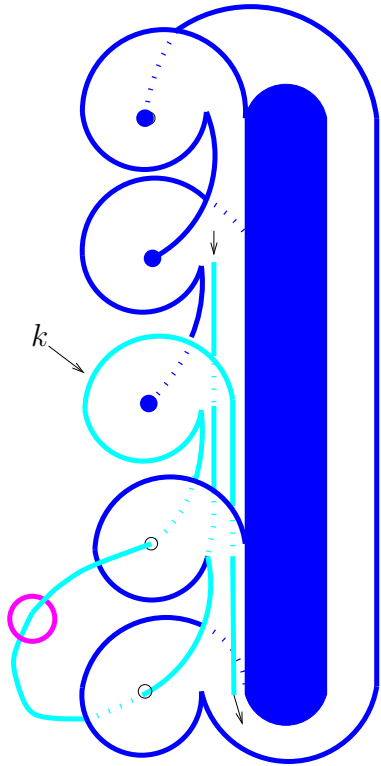}  
		}
	\end{center}
	\caption{The attractor $k$ is inside the light gray tube, 
	the repeller is in the dark gray tube while the dark 
	region is a ball in the basin of attraction of $k$}
	\label{anyk.eps}
\end{figure}

\begin{pf}
The template $U$ shown in Figure~\ref{univ.eps}(a) is known to contain
all knots as periodic orbits \cite{Ghr97}. 
Thus we can suppose  $k$ has been realized as an orbit in $U$.
We shall work with a 
variation of $U$ shown in Figure~\ref{univ.eps}(b) and denoted $V$.
It has five ``Lorenz ears''. Notice however that the middle ear does
not stretch all the way across; it is to extend only as far as an
outer most arc of $k$. (Technically $V$ is not a
template, but it is still a branched manifold with a semi-flow).

Now consider the rather odd looking object in Figure~\ref{anyk.eps}.
The dark gray circle represents the tubular neighborhood of a
repeller. 
The the light gray tube has the same knot type as $k$ (though only
a portion of it is shown); we have only added an 
extra loop in an outermost strand of $k$. The dark region is
a topological ball which meets the light gray tube at a single disk near
the cusp of the fourth Lorenz ear.
Their union is a solid torus $A$. The branched manifold $V$ 
has been cut open along $k$ and is now a true template $T$ (compare
with the proof of Theorem A.3.3 in \cite{GHS}). The boundary of
$T$ is in the boundary of $A$. We thicken up $T$ to get $TT$.
Now we can regard $TT$ as a neighborhood of a saddle set. Its
exit set is attached to $A$ as required. 
From the figure we can see that $A \cup TT$ is a solid unknotted
torus. Thus, we can use a meridian of $k$ as
a repeller and build up the desired flow.
\QED \end{pf}

% @@@@@@@@@@@@@@@@@@@@@@@@@@@@ FIGURE @@@@@@@@@@@@@@@@@@@@@@@@@@@@@@
% @@@@@@@@@@@@@@@@@@@@@@@@@@@@@@@@@@@@@@@@@@@@@@@@@@@@@@@@@@@@@@@@@@@

% @@@@@@@@@@@@@@@@@@@@@@@@@@@@ FIGURE @@@@@@@@@@@@@@@@@@@@@@@@@@@@@@
% @@@@@@@@@@@@@@@@@@@@@@@@@@@@@@@@@@@@@@@@@@@@@@@@@@@@@@@@@@@@@@@@@@@

\begin{cor} \label{cor_alex}
The Alexander polynomial multiplies under connected sums of knots.
\end{cor}

\begin{pf}
The claim is that given knots $k_1$ and $k_2$ then $\Delta_{k_1}
\cdot \Delta_{k_2} = \Delta_{k_1 \# k_2}$.  
By Theorem~\ref{thm_anyk} there exist Smale flows for $k_1$ and 
$k_2$ that satisfy the hypotheses of Theorem~\ref{thm_sums}.
We use conclusion (a) of Theorem~\ref{thm_sums} to construct a
Smale flow with attractor $k_1 \# k_2$ and apply 
\cite[Theorem 4.1]{Fra81} noting that there are only two saddle
sets and hence only two linking matrices needed in the formula
given in \cite[Theorem 4.1]{Fra81}. 
\QED \end{pf}

\begin{rem}[Concluding remarks]
We have in these last two sections given a variety of tools for
building new Smale flows from old ones. Many other such results 
could be stated. But, we are nowhere close to developing a calculus
of Smale flows along the lines of what Wada has done for Morse-Smale
flows. Indeed we don't even know if there are any restrictions on
the link type of $a \cup r$ for simple Smale flows.
\end{rem}

%%%%%%%%%%%%%%%%%%%%%%%%%%%%%%%%%%%
%                                 %
%           BIBLIOGRAPHY          %
%                                 %
%%%%%%%%%%%%%%%%%%%%%%%%%%%%%%%%%%%

%\pagebreak
%\bibliographystyle{plain}
%\bibliography{../BIB/master}

\end{document}